\documentclass[a4paper,12pt]{article}
\usepackage{graphicx}
\usepackage{hyperref}
\usepackage{amsmath}
\usepackage{amssymb}
\usepackage{amsfonts}
\usepackage{url}
\usepackage{fourier}
\def\bc{\begin{center}}
\def\ec{\end{center}}
\def\bi{\begin{itemize}}
\def\ei{\end{itemize}}
\def\be{\begin{equation}}
\def\ee{\end{equation}}

\newcommand{\ie}{{\em i.e.}}  
  
\newcommand{\half}{\frac{1}{2}}  
\begin{document}

\title{Yet another elementary proof of Morley's theorem} 
\author{
St\'ephane Peign\'e \\ 
SUBATECH UMR 6457 (IMT Atlantique, Universit\'e de Nantes, IN2P3/CNRS) \\ 4 rue Alfred Kastler, 44307 Nantes, France
}

\maketitle
\begin{abstract}
I present a simple, elementary proof of Morley's theorem, highlighting the naturalness of this theorem.
\end{abstract}

\section{Introduction}

\begin{figure}[t]
\centering
\includegraphics[width=6cm]{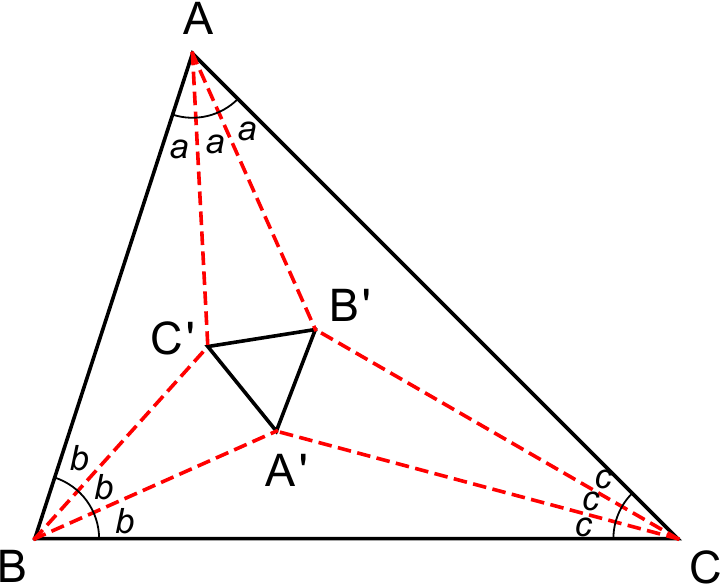}
\caption{Morley's triangle ${\cal M}({\rm ABC})= \rm{A'B'C'}$ of a triangle ABC of angles $3a$, $3b$, $3c$ (with $a+b+c= \frac{\pi}{3}$). According to Morley's theorem, $\rm{A'B'C'}$ is always equilateral.}
\label{fig:morley}
\end{figure}

Morley's theorem states that the trisectors of any (non-degenerate) ABC triangle, taken two by two adjacent to the same side of the triangle, intersect at the vertices of an equilateral triangle, called Morley's triangle ${\cal M}({\rm ABC})$ of the ABC triangle (see Fig.~\ref{fig:morley}).

Among the numerous existing proofs of this theorem~\cite{Weisstein}, we can quote the elementary proof of M.~T.~Naraniengar~\cite{Naraniengar} described for example in the book of Coxeter and Greitzer~\cite{Coxeter}, or more recent proofs, such as that of C.~Frasnay using direct similarities~\cite{Frasnay}, or that of A.~Connes using complex numbers~\cite{Connes}.~The `extraordinary' or `miraculous' nature of this theorem has often been noted~\cite{Coxeter,Lange,Newman}. 

The construction of ${\cal M}({\rm ABC})$ starting from ABC depends only on angles and does not use any notion of metrics. Therefore, if ${\cal M}({\rm ABC})$ is equilateral, ${\cal M}({s}{\rm(ABC)})$ is a fortiori equilateral for any triangle $s{\rm(ABC)}$ obtained from ABC by a direct similarity $s$. Morley's theorem can thus be reformulated by restricting oneself, for each {\it shape} of the ABC triangle, to only one {\it size} of this triangle. The latter can be freely chosen and adjusted so that ${\cal M}({\rm ABC})$ is always superimposable on the same equilateral triangle. This leads to the following formulation (sometimes used implicitly) of Morley's theorem:

\vspace{3mm}
{\it For a given equilateral triangle ${\rm A'B'C'}$ and for any choice $(a,b,c)$ of positive (non-zero) real numbers verifying $a+b+c= \frac{\pi}{3}$, we can construct a triangle ABC of angles $3a$, $3b$, $3c$ such that ${\cal M}({\rm ABC})= {\rm A'B'C'}$.} 
\vspace{3mm}

In this paper, I present a proof of Morley's theorem stated in this form. One objective of the proposed demonstration (apparently not present in the literature) is to make this theorem natural. 

\section{A natural proof of Morley's theorem}
\label{sec:Morley} 

\begin{figure}[t]
\centering
\includegraphics[width=7cm]{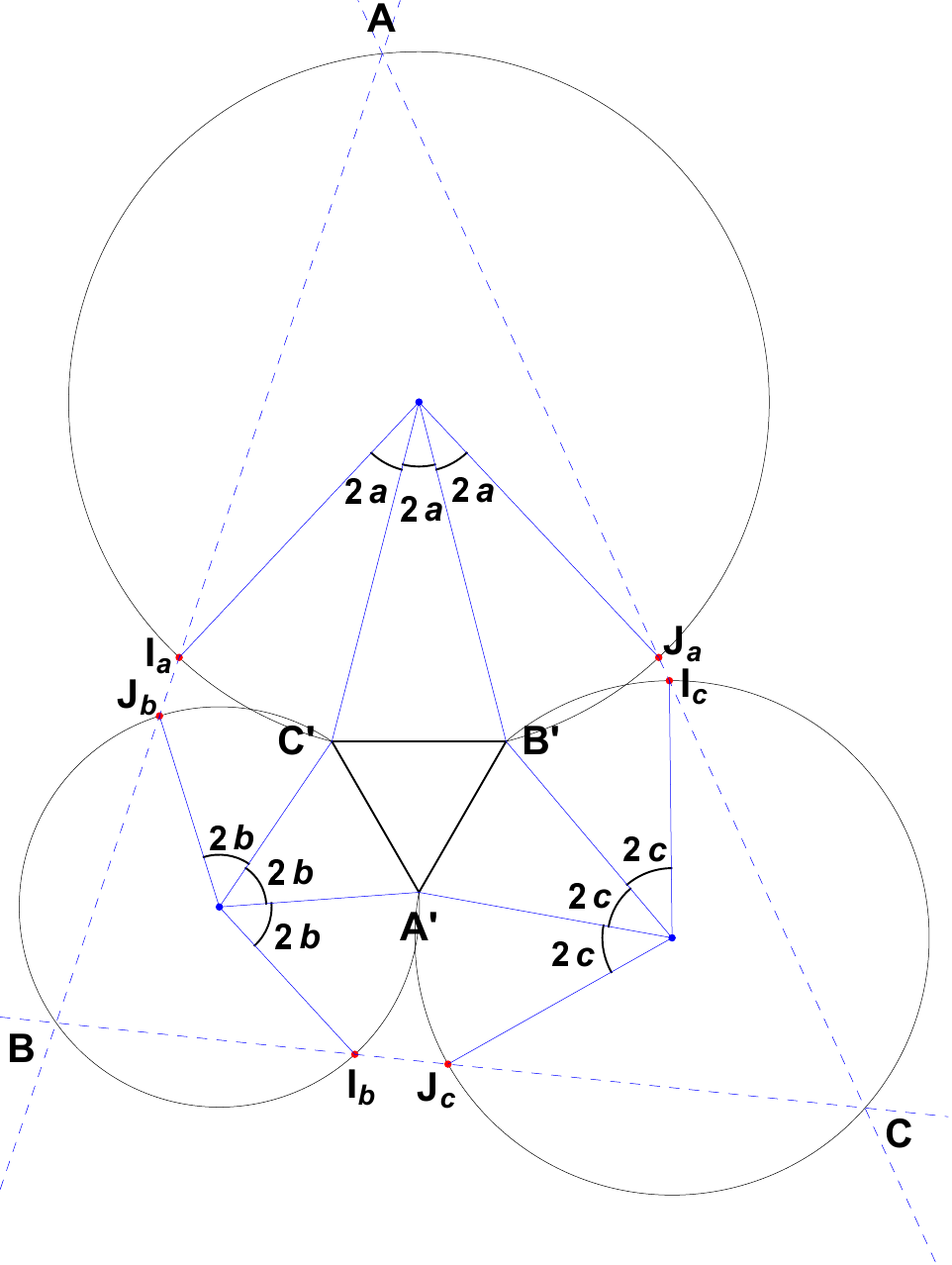}
\vspace{-3mm}
\caption{Construction, starting from a given equilateral triangle ${\rm A'B'C'}$, of an ABC triangle of angles $3a$, $3b$, $3c$ satisfying ${\cal M}({\rm ABC})= \rm{A'B'C'}$.}
\label{fig:construction}
\end{figure}

Let us thus start from an equilateral triangle ${\rm A'B'C'}$ and given strictly positive numbers $a$, $b$, $c$ such that $a+b+c= \frac{\pi}{3}$, and let us show that it is possible to construct an ABC triangle of angles $3a$, $3b$, $3c$ having ${\rm A'B'C'}$ for Morley's triangle. The proof is in two steps.
\bi
\item[] (i) We first show that the points A, B and C, if they exist, can only be located at positions determined by a quasi-automatic construction (thus requiring no tricks that might make this part of the proof look somewhat miraculous). This construction is as follows:
\ei
\vspace{-2mm}
Since we must have $\widehat{\rm C'AB'}= a$ (see Fig.~\ref{fig:morley}), the point A must belong to the circular arc $\wideparen{{\rm B'}{\rm C'}}$ such that the chord $\rm{C'B'}$ subtends a central angle equal to $2a$. The angle $a$ being given, this arc (a fortiori located on the opposite side of ${\rm A'}$ with respect to the line $\rm{(C'B')}$) is easily constructed (see Fig.~\ref{fig:construction}). We must also have $\widehat{\rm BAC'}=\widehat{\rm B'AC}= a$ (Fig.~\ref{fig:morley}). The lines $\rm{(AB)}$ and $\rm{(AC)}$ must therefore intersect $\wideparen{{\rm B'}{\rm C'}}$ respectively at the points ${\rm I}_a$ and ${\rm J}_a$ forming the inscribed angles $\widehat{{\rm I}_a{\rm AC'}}= a$ and $\widehat{{\rm B'A}{\rm J}_a}= a$, \ie\ such that the chords ${\rm I}_a\rm{C'}$ and $\rm{B'}{\rm J}_a$ each subtend a central angle equal to $2a$. The points ${\rm I}_a$ and ${\rm J}_a$ are thus easily constructed knowing $a$ (Fig.~\ref{fig:construction}). The same reasoning obtained by circular permutation of the triplets $({\rm A}, {\rm B}, {\rm C})$, $({\rm A'}, {\rm B'}, {\rm C'})$, and $(a, b, c)$ allows one to construct the arc $\wideparen{{\rm C'}{\rm A'}}$ and the points ${\rm I}_b$ and ${\rm J}_b$, then the arc $\wideparen{{\rm A'}{\rm B'}}$ and the points ${\rm I}_c$ and ${\rm J}_c$, and to conclude that the lines $\rm{(AB)}$, $\rm{(BC)}$, and $\rm{(CA)}$ must coincide respectively with the lines $({\rm I}_a{\rm J}_b)$, $({\rm I}_b{\rm J}_c)$, and $({\rm I}_c{\rm J}_a)$. Therefore, the points A, B, C are {\it necessarily} located at the intersections of these lines.
\bi 
\item[] (ii) We now show that this condition is {\it sufficient}, \ie\ that A, B, C {\it defined} by ${\rm A} \equiv ({\rm I}_a {\rm J}_b) \cap ({\rm I}_c{\rm J}_a)$, ${\rm B} \equiv ({\rm I}_b{\rm J}_c) \cap ({\rm I}_a{\rm J}_b)$, and ${\rm C} \equiv ({\rm I}_c  {\rm J}_a) \cap ({\rm I}_b{\rm J}_c)$, are such that ${\cal M}({\rm ABC})= \rm{A'B'C'}$. For this it is sufficient to show that $\widehat{{\rm I}_a {\rm A} {\rm J}_a} = 3 a$, which simply follows from the previous construction:
\ei
\vspace{-2mm}
The sum of the angles around the point B' being equal to $2 \pi$, we find $\widehat{{\rm I}_c {\rm B'} {\rm J}_a} = \frac{\pi}{3} -2b$. The triangle ${\rm I}_c {\rm B'} {\rm J}_a$ being isosceles, we infer $\widehat{{\rm A}{\rm J}_a{\rm B'}} = \half (\pi + \widehat{{\rm I}_c {\rm B'} {\rm J}_a}) = \frac{2\pi}{3} -b$. By symmetry, we have $\widehat{{\rm A}{\rm I}_a{\rm C'}} = \frac{2\pi}{3} -c$. Finally, the sum of the angles of the pentagon ${\rm A}{\rm I}_a{\rm C'}{\rm B'}{\rm J}_a$ equalling $3\pi$, we get $\widehat{{\rm I}_a {\rm A} {\rm J}_a} = 3a$. Thus, the point A belongs to the arc $\wideparen{{\rm J}_a {\rm I}_a}$ whose chord ${\rm I}_a {\rm J}_a$ subtends a central angle equal to $6a$, which is a portion of the arc $\wideparen{{\rm B'}{\rm C'}}$ constructed in step (i). The lines $\rm{(AC')}$ and $\rm{(AB')}$ are thus the trisectors of $\widehat{{\rm B}{\rm A}{\rm C}}$. Similar assertions for the angles at B and C of ABC are obtained by circular permutation, demonstrating Morley's theorem.

\vspace{3mm}
Note that the central point of the proof, $\widehat{{\rm I}_a {\rm A} {\rm J}_a} = 3a$, can also be deduced from the symmetry properties of the construction, without going through the above explicit (admittedly trivial) calculation of $\widehat{{\rm I}_a {\rm A} {\rm J}_a}$. Indeed, by construction $\widehat{{\rm I}_a {\rm A} {\rm J}_a}$ must be a {\it linear} function of $a, b, c$, {\it symmetric} in the exchange $b \leftrightarrow c$. Therefore, $\widehat{{\rm I}_a {\rm A} {\rm J}_a}$ is a linear function of $a$ and $b+c$, and since $a+b+c= \frac{\pi}{3}$, a linear function of $a$ only, $\widehat{{\rm I}_a {\rm A} {\rm J}_a} = f(a) = \alpha a + \beta$. By symmetry we must have $\widehat{{\rm I}_b {\rm B} {\rm J}_b} = f(b)$ and $\widehat{{\rm I}_c {\rm C} {\rm J}_c} = f(c)$, with the same function $f$. This function is easily determined by a mere observation of the construction. The constant $\beta$ is obtained by considering the limit of $f(a)$ when $a \to 0$. Since $f(a)$ depends only on $a$, we can set  $b=c$ to determine this limit, and thus $2b = \frac{\pi}{3} - a \to \frac{\pi}{3}$. In this limit, we see on Fig.~\ref{fig:construction} that the points ${\rm I}_a$ and ${\rm J}_b$ tend to the symmetrical point S of ${\rm B'}$ with respect to ${\rm C'}$, while remaining on the circle of centre ${\rm C'}$ and radius ${\rm C'}{\rm B'}$. The line $({\rm I}_a{\rm J}_b)$ thus tends towards the tangent to this circle at point S, which is perpendicular to $({\rm C'}{\rm B'})$. Similarly, $({\rm I}_c{\rm J}_a)$ becomes perpendicular to $({\rm C'}{\rm B'})$. We infer that $\widehat{{\rm I}_a {\rm A} {\rm J}_a} \to 0$ when $a \to 0$, hence $\beta =0$. Finally, $f(a)+f(b)+f(c)= \pi$ implies $\alpha = 3$, and thus $\widehat{{\rm I}_a {\rm A} {\rm J}_a} = 3 a$. 

\vspace{3mm}
By way of conclusion, let us mention the proof of Morley's theorem by J.~Conway~\cite{Conway}, which consists in giving oneself {\it at the start} six triangles with well-chosen shapes (which can be checked from Ref.~\cite{Conway} to coincide with those of the triangles ${\rm ABC'}$, ${\rm BCA'}$, ${\rm CAB'}$, ${\rm AC'B'}$, ${\rm BA'C'}$, ${\rm CB'A'}$ of Fig.~\ref{fig:construction}), then adjusting their sizes, and finally showing that the obtained triangles fit perfectly around the equilateral triangle ${\rm A'B'C'}$ to form the sought ABC triangle. This provides an elegant and concise proof, which however requires some substantial amount of intuition. The proof presented here explains where these triangles come from, and why the constructed triangle ABC satisfies ${\cal M}({\rm ABC})= \rm{A'B'C'}$ from simple symmetry arguments. We hope it will help demystify `Morley's miracle'.

\providecommand{\href}[2]{#2}\begingroup\raggedright
\endgroup

\end{document}